# Baker's dozen digits of two sums involving reciprocal products of an integer and its greatest prime factor

Tengiz O. Gogoberidze

**Abstract.** Two sums $\sum_{n=1}^{\infty} \frac{d(n)}{n \cdot G(n)}$ and $\sum_{n=1}^{\infty} \frac{1}{n \cdot G(n)}$ over the inverse of the product of an integer $n$ and its greatest prime factor $G(n)$, are computed to first 13 decimal digits. These sums converge, but converge very slowly. They are transformed into sums involving Mertens' prime product $\prod_{p \leq x} \frac{p}{p-1}$ with the remainder term which are estimated by means of Chebyshev's $\vartheta$-function.

## 1. Introduction

Let $G(n)$ denote the greatest prime factor of an integer $n \geq 2$, $G(1) = 1$ [11, A006530], and let $d(n) = \sigma_0(n)$ be the number of divisors of $n$ [11, A000005].

Divergent sums involving reciprocals of $G(n)$ have been investigated in many works and some important estimates concerning asymptotic behavior of these sums were obtained [12], [13].

The purpose of this paper is to calculate some first decimal digits of two convergent sums:

$$Sa \equiv \sum_{n=1}^{\infty} \frac{d(n)}{n \cdot G(n)} \tag{1.1}$$

and

$$Sb \equiv \sum_{n=1}^{\infty} \frac{1}{n \cdot G(n)}. \tag{1.2}$$

Thinking "naive" one can predict behavior of the sum (1.1) based on mean values $\langle G(n) \rangle = \frac{1}{n} \cdot \sum_{k=1}^{n} G(k) \sim \frac{\pi^2}{12} \cdot \frac{n}{\ln n}$ [1], $\langle d(n) \rangle = \frac{1}{n} \cdot \sum_{k=1}^{n} d(k) \sim \ln n$ [3]. The sum $\sum_{n=1}^{\infty} \left(\frac{\ln n}{n}\right)^2 = \zeta''(2) = 1.98928\ldots$ [7, Table 2] is well-known, however, $Sa$ convergence is really slow. Table 1 shows partial sums $Sa_n$ for some powers of 10.

Table 1. Partial sums $Sa_n$ for some powers of 10

| $n$ | $Sa_n$ |
|---|---|
| $10^2$ | 4.816507 … |
| $10^3$ | 6.149878 … |
| $10^4$ | 6.961411 … |
| $10^5$ | 7.4338008 … |
| $10^6$ | 7.7094259 … |
| $10^7$ | 7.870104 … |

These partial sums can be approximated by the fitting curve $Sa_n = 8.115653 - 9.327239 \cdot n^{-0.226343}$. So to obtain first decimal digit adding up $2.7 \cdot 10^8$ terms may be required.

As to (1.2), it converges somewhat more rapidly. For 3 true decimal digits $3 \cdot 10^6$ terms of series are quite enough: $Sb_{3 \cdot 10^6} = \mathbf{2.25}05789 \ldots$ Therefore, it is convenient to start counting with second sum.

## 2. Transformation of the sums and some preliminary results

Since $\prod_{p \le x} \frac{p}{p-1} = \sum_{G(n) \le x} \frac{1}{n}$ [14, (1.4)], one can write

$$Sb = \sum_{n=1}^{\infty} \frac{1}{n \cdot G(n)} = 1 + \sum_p \left[\frac{1}{p} \cdot \sum_{G(n)=p} \frac{1}{n}\right] = 1 + \sum_p \left[\frac{1}{p} \cdot \left\{\sum_{G(n) \le p} \frac{1}{n} - \sum_{G(n) < p} \frac{1}{n}\right\}\right]$$

$$= 1 + \sum_p \left[\frac{1}{p} \cdot \prod_{p' \le p} \frac{p'}{p'-1}\left(1 - \frac{p-1}{p}\right)\right] = 1 + \sum_p \left[\frac{1}{p^2} \cdot \prod_{p' \le p} \frac{p'}{p'-1}\right].$$

Here a product $\prod_{p' \le p} \frac{p'}{p'-1}$ is taken over all primes $p' \le p$.

In 1874, F. Mertens showed that $\prod_{p \le x} \frac{p}{p-1} \sim e^\gamma \ln x$ [2, p.53], where $\gamma = 0.57721 \ldots$ is the Euler-Mascheroni constant [11, A001620].

Owing to Mertens' product asymptotics, transformation process can be continued as follows:

$$Sb = 1 + e^\gamma \sum_p \frac{\ln p}{p^2} + \sum_{p \le x} \left[\frac{1}{p^2}\left\{\prod_{p' \le p} \frac{p'}{p'-1} - e^\gamma \ln p\right\}\right] + \sum_{p > x} \left[\frac{1}{p^2}\left\{\prod_{p' \le p} \frac{p'}{p'-1} - e^\gamma \ln p\right\}\right]$$

$$Sb = Cb + \sum_{p \le x} \left[\frac{1}{p^2}\left\{\prod_{p' \le p} \frac{p'}{p'-1} - e^\gamma \ln p\right\}\right] + Rb(x), \tag{2.1}$$

where $Cb = 1 + e^\gamma \sum_p \frac{\ln p}{p^2} = 1 - e^\gamma P'(2) = 1.878230974452894488025 3707 \ldots$ is related to prime zeta function derivative $P'(2) = -\sum_p \frac{\ln p}{p^2} = -0.4930911093687644621 97826205 \ldots$ [6, Table 3], [11, A136271], and $Rb(x) = \sum_{p > x} \left[\frac{1}{p^2}\left\{\prod_{p' \le p} \frac{p'}{p'-1} - e^\gamma \ln p\right\}\right]$ is remainder term.

Analogous transformation of the $Sa$ sum is shown below:

$$Sa = \sum_{n=1}^{\infty} \frac{d(n)}{n \cdot G(n)} = 1 + \sum_p \left[\frac{1}{p} \cdot \left(\sum_{G(n)=p} \frac{1}{n}\right)^2\right] = 1 + \sum_p \left[\frac{1}{p} \cdot \left\{\sum_{G(n) \le p} \frac{1}{n} - \sum_{G(n) < p} \frac{1}{n}\right\}^2\right]$$

$$= 1 + \sum_p \left[\frac{1}{p} \cdot \left(\prod_{p' \le p} \frac{p'}{p'-1}\right)^2 \cdot \left(1 - \left(\frac{p-1}{p}\right)^2\right)\right] = 1 + \sum_p \left[\frac{2 - \frac{1}{p}}{p^2} \cdot \left(\prod_{p' \le p} \frac{p'}{p'-1}\right)^2\right]$$

$$= 1 + e^{2\gamma} \sum_p \frac{\ln^2 p}{p^2}\left(2 - \frac{1}{p}\right) + \sum_p \left[\frac{2 - \frac{1}{p}}{p^2}\left\{\left(\prod_{p' \le p} \frac{p'}{p'-1}\right)^2 - e^{2\gamma} \ln^2 p\right\}\right].$$

Then

$$Sa = Ca + \sum_{p \le x} \left[\frac{2 - \frac{1}{p}}{p^2}\left\{\left(\prod_{p' \le p} \frac{p'}{p'-1}\right)^2 - e^{2\gamma} \ln^2 p\right\}\right] + Ra(x), \tag{2.2}$$

where $Ca = 1 + e^{2\gamma}\left[2 \sum_p \frac{\ln^2 p}{p^2} - \sum_p \frac{\ln^2 p}{p^3}\right] = 1 + e^{2\gamma}[2 \cdot P''(2) - P''(3)]$ is related to prime zeta function $2^{nd}$ derivative $P''(s)$, and $Ra(x) = \sum_{p > x} \left[\frac{2 - \frac{1}{p}}{p^2}\left\{\left(\prod_{p' \le p} \frac{p'}{p'-1}\right)^2 - e^{2\gamma} \ln^2 p\right\}\right]$ is remainder term.

## 3. Statement of main results

**Theorem 1.**

$$Sb = Cb + \sum_{p \leq x} \left[ \frac{1}{p^2} \left\{ \prod_{p' \leq p} \frac{p'}{p'-1} - e^\gamma \ln p \right\} \right] + Rb(x),$$

where $Cb = 1 - e^\gamma P'(2) = 1.878230974452894488025 3707\ldots$, and $-\frac{0.034}{x \ln^3 x} < Rb(x) < \frac{0.1}{x \ln^3 x}$ for every $x \geq 51841229$.

**Corollary 1.1**

$$Sb = 2.254435359519\ldots$$

*Proof.* For $x = 86028161$ as a result of PC calculation we have $Sb = 2.254435359519071\ldots + Rb(86028161)$, $-0.6 \cdot 10^{-13} < Rb(86028161) < 2 \cdot 10^{-13}$.

**Theorem 2.**

$$Sa = Ca + \sum_{p \leq x} \left[ \frac{2 - \frac{1}{p}}{p^2} \left\{ \left( \prod_{p' \leq p} \frac{p'}{p'-1} \right)^2 - e^{2\gamma} \ln^2 p \right\} \right] + Ra(x),$$

where $Ca = 1 + e^{2\gamma}[2 \cdot P''(2) - P''(3)] = 5.229250296762544252 9818603\ldots$, and $-\frac{0.24}{x \ln^2 x} < Ra(x) < \frac{0.72}{x \ln^2 x}$ for every $x \geq 51841229$.

**Corollary 2.1**

$$Sa = 8.115653111459\ldots$$

*Proof.* For $x = 2576983867$ as a result of PC calculation we have $Sa = 8.115653111459203\ldots + Ra(2576983867)$, $-2 \cdot 10^{-13} < Ra(2576983867) < 6 \cdot 10^{-13}$.

## 4. Estimates for Mertens' prime product

A lot of modern estimates for some functions over primes, including Mertens' product, are based on the method investigated by J. B. Rosser and L. Schoenfeld. They connect the sum of reciprocal primes not exceeding $x$ with Chebyshev's $\vartheta$-function in the following way [4, (4.20)]:

$$\sum_{p \leq x} \frac{1}{p} - \ln \ln x - B = \frac{\vartheta(x) - x}{x} - \int_x^\infty \frac{\vartheta(y) - y}{y} \cdot \frac{1 + \ln y}{y \ln^2 y} dy,$$

where $B = \gamma + \sum_p \left\{ \ln\left(1 - \frac{1}{p}\right) + \frac{1}{p} \right\} = 0.2614972128\ldots$ – Mertens' constant [11, A077761].

Using this expression and the form of $\vartheta$-function bounds like $\left|\frac{\vartheta(x)-x}{x}\right| < \frac{\eta_k}{\ln^k x}$, P. Dusart derived the inequality

$$\left| \sum_{p \leq x} \frac{1}{p} - \ln \ln x - B \right| \leq \frac{\eta_k}{k \cdot \ln^k x} + \eta_k \frac{1 + \frac{1}{k+1}}{\ln^{k+1} x}, \qquad [9, (5.7)]$$

which helped him to improve "classical" Rosser-Schoenfeld estimate [4, Theorem 8] for Mertens' prime product

$$e^\gamma \ln(x) \cdot \left\{ 1 - \frac{0.2}{\ln^3 x} \right\} < \prod_{p \leq x} \frac{p}{p-1} < e^\gamma \ln(x) \cdot \left\{ 1 + \frac{0.2}{\ln^3 x} \right\}, \qquad (4.1)$$

for every $x \geq 2278382$ [9, Theorem 5.9].

Also, more recently C. Axler [10, Prop. 6.1] showed that

$$\frac{e^{-\gamma}}{\ln(x)} \cdot \left\{1 - \frac{0.05}{\ln^3 x} - \frac{3}{16\ln^4 x}\right\} < \prod_{p \leq x} \frac{p-1}{p} < \frac{e^{-\gamma}}{\ln(x)} \cdot \left\{1 + \frac{0.07}{\ln^3 x}\right\},$$

where the left-hand side inequality is valid for every $x \geq 46909074$ and the right-hand side inequality is valid for every $x > 1$.

Finally, he proposed even sharper estimate for $x \geq 1797126630797 = p_{64707865143}$, but summarizing over so many primes would take a long computer time for our sums calculation.

To obtain a more usable estimate for Mertens' product let us take advantage of an "asymmetrical" form of Chebyshev's $\vartheta$-function bounds

$$-\frac{\eta_k^-}{\ln^k x} < \frac{\vartheta(x) - x}{x} < \frac{\eta_k^+}{\ln^k x}. \tag{4.2}$$

So, according to Broadbent et al. [8, Col. 11.1], $\eta_3^- = 0.15$ for $x \geq 19035709163$; $\eta_3^+ = 0.024334$ for $x > 1$. Then inequality [9, (5.7)] can be rewritten like

$$-\left\{\frac{\eta_k^+}{k \cdot \ln^k x} + \frac{\frac{\eta_k^+}{k+1} + \eta_k^-}{\ln^{k+1} x}\right\} < \sum_{p \leq x} \frac{1}{p} - \ln\ln x - B < \frac{\eta_k^-}{k \cdot \ln^k x} + \frac{\frac{\eta_k^-}{k+1} + \eta_k^+}{\ln^{k+1} x} \tag{4.3}$$

Following Dusart's way, we obtain several successive inequalities

$$-\Sigma^-(x) + S(x) < -\sum_{p \leq x} \ln\left(1 - \frac{1}{p}\right) - \ln\ln x - \gamma < \Sigma^+(x) + S(x),$$

where $\Sigma^-(x)$ and $\Sigma^+(x)$ denote absolute values of left- and right-hand sides of (4.3), respectively, and $S(x) = \sum_{p > x}\left\{\ln\left(1 - \frac{1}{p}\right) + \frac{1}{p}\right\}$ is such that $-\frac{1.02}{(x-1)\ln x} < S(x) < 0$ by Rosser and Schoenfeld.

$$e^\gamma \ln(x) \cdot \exp\left\{-\Sigma^-(x) - \frac{1.02}{(x-1)\ln x}\right\} < \prod_{p \leq x} \frac{p}{p-1} < e^\gamma \ln(x) \cdot \exp\{\Sigma^+(x)\}$$

Or, for $x \geq x_0$

$$e^\gamma \ln(x) \cdot \left\{1 - \Sigma^-(x) - \frac{1.02}{(x-1)\ln x}\right\} < \prod_{p \leq x} \frac{p}{p-1} < e^\gamma \ln(x) \cdot \left\{1 + \frac{e^{\Sigma^+(x_0)} - 1}{\Sigma^+(x_0)} \cdot \Sigma^+(x)\right\}, \tag{4.4}$$

because of $\left(\frac{e^z - 1}{z}\right)' > 0$ and $\Sigma^+(x) \leq \Sigma^+(x_0)$.

Using Dusart's result (4.1) and estimate [8, Col. 11.1], we continue the right-hand side of inequality (4.4) and find

$$\prod_{p \leq x} \frac{p}{p-1} < e^\gamma \ln(x) \cdot \left\{1 + 1.0000319 \cdot \left[\frac{0.05}{\ln^3 x} + \frac{0.061834}{\ln^4 x}\right]\right\} < e^\gamma \ln(x) \cdot \left\{1 + \frac{0.0561}{\ln^3 x}\right\}.$$

Analogously we obtain on the left-hand side of (4.4)

$$\prod_{p \leq x} \frac{p}{p-1} > e^\gamma \ln(x) \cdot \left\{1 - \left[\frac{0.024334}{3 \cdot \ln^3 x} + \frac{0.1560835}{\ln^4 x}\right] - \frac{1.02}{(x-1)\ln x}\right\} > e^\gamma \ln(x) \cdot \left\{1 - \frac{0.0189}{\ln^3 x}\right\}.$$

Both inequalities holds at least for $x \geq 19035709163$.

As a result, we get

$$e^\gamma \ln(x) \cdot \left\{1 - \frac{0.0189}{\ln^3 x}\right\} < \prod_{p \leq x} \frac{p}{p-1} < e^\gamma \ln(x) \cdot \left\{1 + \frac{0.0561}{\ln^3 x}\right\} \tag{4.5}$$

for $x \geq 51841229$. To confirm that inequality (4.5) holds for $51841229 \leq x < 19035709163$, computer check has been used.

## 5. Useful lemma

Applying (4.5) and $Rb(x)$ definition by (2.1), we have

$$-\frac{0.0189 \cdot e^\gamma}{\ln^2 x} < \prod_{p \leq x} \frac{p}{p-1} - e^\gamma \ln(x) < \frac{0.0561 \cdot e^\gamma}{\ln^2 x},$$

$$-\sum_{p>x} \frac{0.034}{p^2 \ln^2 p} < Rb(x) < \sum_{p>x} \frac{0.1}{p^2 \ln^2 p} \tag{5.1}$$

and

$$-\frac{0.034}{x \ln^3 x} < Rb(x) < \frac{0.1}{x \ln^3 x} \tag{5.2}$$

To substantiate the transition from (5.1) to (5.2), it suffices to prove following lemma.

**Lemma 5.1**

$$\sum_{p>x} \frac{1}{p^2 \ln^m p} < \frac{1}{x \ln^{m+1} x} \cdot \left\{ 1 - \frac{m+1}{\ln x} + \left(\frac{m+3/2}{\ln x}\right)^2 \right\},$$

for every $x \geq 5 \cdot 10^7; m > 0$.

*Proof.* The expression below proposed by Rosser and Schoenfeld allows finding some sums over primes.

$$\sum_{p \leq x} f(p) = \frac{f(x) \cdot \vartheta(x)}{\ln x} - \int_2^x \vartheta(y) \left(\frac{f(y)}{\ln y}\right)' dy, \qquad [4.(4.13)]$$

Denoting $\frac{1}{p^2 \ln^m p}$ by $f(p)$, we get

$$\sum_{p>x} f(p) = -\frac{\vartheta(x)}{x^2 \ln^{m+1} x} + \int_2^x \vartheta(y) \left(\frac{1}{y^2 \ln^{m+1} y}\right)' dy = -\frac{\vartheta(x) - x}{x^2 \ln^{m+1} x}$$

$$+ \left\{ -\frac{1}{x \ln^{m+1} x} + 2 \int_2^x \frac{dy}{y^2 \ln^{m+1} y} + (m+1) \int_x^\infty \frac{dy}{y^2 \ln^{m+2} y} \right\} + \int_x^\infty \frac{\vartheta(y) - y}{y} \cdot \frac{2 + \frac{m+1}{\ln y}}{y^2 \ln^{m+1} y} dy$$

Since $\left|\frac{\vartheta(x)-x}{x}\right| < \frac{\eta_2}{\ln^2 x}$, integration by parts gives

$$\sum_{p>x} f(p) < \frac{1}{x \ln^{m+1} x} - (m+1) \int_x^\infty \frac{dy}{y^2 \ln^{m+2} y} + \eta_2 \left\{ \frac{1}{x \ln^{m+3} x} + 2 \int_2^x \frac{dy}{y^2 \ln^{m+3} y} + \int_x^\infty \frac{m+1}{y^2 \ln^{m+4} y} dy \right\},$$

and once more

$$\sum_{p>x} f(p) < \frac{1}{x \ln^{m+1} x} \left\{ 1 - \frac{m+1}{\ln x} + \frac{(m+1)(m+2)}{\ln^2 x} \right\} - (m+1)(m+2)(m+3) \int_x^\infty \frac{dy}{y^2 \ln^{m+4} y}$$

$$+ \eta_2 \left\{ \frac{3}{x \ln^{m+3} x} - (m+5) \int_x^\infty \frac{dy}{y^2 \ln^{m+4} y} \right\} < \frac{1}{x \ln^{m+1} x} \left\{ 1 - \frac{m+1}{\ln x} + \frac{(m+1)(m+2) + 3\eta_2}{\ln^2 x} \right\}.$$

According to S. Broadbent et al. [8, Col. 11.1], $\eta_2 = 0.068701$ for $x \geq 5 \cdot 10^7$ [8, Table 15], hence

$$\sum_{p>x} f(p) < \frac{1}{x \ln^{m+1} x} \left\{ 1 - \frac{m+1}{\ln x} + \frac{(m+1)(m+2) + 0.207}{\ln^2 x} \right\}.$$

Q.E.D. **Theorem 1**, therefore, has also been proven.

## 6. Prime zeta function 2$^{nd}$ derivative

To investigate prime zeta function $P(s) \equiv \sum_p^\infty \frac{1}{p^s}$ and its derivatives, it's frequently convenient to follow procedure [5] based on Möbius inversion of Riemann zeta function logarithmic derivative.

Introducing partial function like
$$\zeta_{p>x}(s) \equiv \zeta(s) \cdot \prod_{p \leq x}\left(1 - \frac{1}{p^s}\right),$$
H. Cohen got a helpful expression for prime zeta function
$$P(s) = \sum_{p \leq x} \frac{1}{p^s} + \sum_{k=1}^\infty \frac{\mu(k)}{k} \ln \zeta_{p>x}(ks),$$
where $\mu(k)$ is Möbius function [11, A008683].

By differentiating of these expressions, R. J. Mathar deduced a formula for quick calculating prime zeta function derivative
$$P'(s) = -\sum_{p \leq x} \frac{\ln p}{p^s} + \sum_{k=1}^\infty \mu(k)\left[\frac{\zeta'(ks)}{\zeta(ks)} + \sum_{p \leq x} \frac{\ln p}{p^{ks} - 1}\right]. \qquad [6, (6); (7)]$$

Differentiating again, we can derive an expression for calculating prime zeta function 2$^{nd}$ derivative
$$P''(s) = \sum_{p \leq x} \frac{\ln^2 p}{p^s} + \sum_{k=1}^\infty k \cdot \mu(k) \cdot \left[\frac{\zeta''(ks)}{\zeta(ks)} - \left(\frac{\zeta'(ks)}{\zeta(ks)}\right)^2 - \sum_{p \leq x} \frac{\ln^2 p}{p^{ks} + p^{-ks} - 2}\right]. \qquad (6.1)$$

Table 2 shows values of $P''(s)$ at some integer arguments $s$.

Table 2. Values of $P''(s)$ at some integer arguments $s$

| $s$ | $P''(s)$ |
| --- | --- |
| 2 | 0.7415978549828050030239403275450375 ... |
| 3 | 0.1499805843420868545613264915839837 ... |
| 4 | 0.0515291349877069852843053881647288 ... |
| 5 | 0.0211008945931815396107750284303329 ... |
| 6 | 0.0093659030058415755113780509464443 ... |

## 7. Proof of theorem 2

First, we substitute the $P''(2)$ and $P''(3)$ values into the expression for $Ca$ constant and get the desired outcome
$$Ca = 5.229250296762544252981860349634 ... \qquad (7.1)$$
Then, applying (4.5) and $Ra(x)$ definition by (2.2), we have
$$e^{2\gamma} \ln^2(x) \cdot \left\{1 - \frac{0.0378}{\ln^3 x}\right\} < \left(\prod_{p \leq x} \frac{p}{p-1}\right)^2 < e^{2\gamma} \ln^2(x) \cdot \left\{1 + \frac{0.1123}{\ln^3 x}\right\},$$
$$-\frac{0.12}{\ln x} < \left(\prod_{p \leq x} \frac{p}{p-1}\right)^2 - e^{2\gamma} \ln^2(x) < \frac{0.36}{\ln x}$$

and

$$-\sum_{p>x} \frac{0.24}{p^2 \ln p} < Ra(x) < \sum_{p>x} \frac{0.72}{p^2 \ln p} \qquad (7.2)$$

for every $x \geq 51841229$.

Note that, by using lemma 5.1, we obtain the remainder term estimate

$$-\frac{0.24}{x \ln^2 x} < Ra(x) < \frac{0.72}{x \ln^2 x}. \qquad (7.3)$$

The theorem is proved.